\documentclass{lmcs}
\usepackage[utf8]{inputenc}

\usepackage{color}
\usepackage{hyperref}
\usepackage{float}
\usepackage{amsmath}
\usepackage{amsfonts}
\usepackage{amsthm}
\usepackage{amssymb}
\usepackage{proof}
\usepackage{mathpartir}
\usepackage{mathrsfs}
\usepackage{stmaryrd}
\usepackage{cmll}
\usepackage{enumerate}
\usepackage{url}

\usepackage{listings}

\newcommand {\q}[2]{{\tt q}_{#1, #2}}
\newcommand {\qI}{{\tt q}}

\newcommand{\comp}[5]{{\tt comp}(#1, #2, #3, #4, #5)}

\newcommand{\Ctx}{\mathrm{Ctx}}
\newcommand{\Sub}{\mathrm{Sub}}
\newcommand{\Ty}{\mathrm{Ty}}
\newcommand{\Tm}{\mathrm{Tm}}
\newcommand{\C}{{\mathcal C}}

\newcommand{\T}{{\mathcal T}}

\newcommand{\RawCtx}{{\tt Ctx}}
\newcommand{\RawSub}{{\tt Sub}}
\newcommand{\RawTy}{{\tt Ty}}
\newcommand{\RawTm}{{\tt Tm}}

\newcommand{\inte}[1]{\llbracket #1 \rrbracket}

\newcommand{\Set}[0]{{\bf Set}}

\newcommand{\omitthis}[1]{}

\newcommand{\changenote}[1]{}

\newcommand{\longtext}[1]{}
\newcommand{\shorttext}[1]{}
\newcommand{\commentaway}[1]{}

\definecolor{Red}{rgb}{1,0,0}

\newcommand{\Fam}{\textbf{Fam}}

\newcommand{\cext}{.}

\newcommand{\op}{\text{op}}

\def\N{\mathsf{N}}
\def\U{\mathsf{U}}

\def\app{\mathsf{app}}
\def\Cop{\C^\op}

\def\p{\mathrm{p}}
\def\q{\mathrm{q}}
\newcommand{\tuple}[1]{\langle #1 \rangle}

\newtheorem{remark}{Remark}
\newtheorem{definition}{Definition}

\newtheorem{theorem}{Theorem}
\newcommand{\s}{\mathrm{s}}
\newcommand{\Rec}{\mathrm{R}}
\newcommand{\Ta}{\mathrm{T}}

\def\D{\mathcal{D}}
\def\V{\mathrm{V}}
\def\Cwf{\mathbf{CwF}}
\def\Obj{\mathrm{obj}}
\def\Ctx{\mathrm{Ctx}}
\def\Hom{\mathrm{hom}}
\def\id{\mathrm{id}}
\def\Mon{\mathrm{M}}
\def\idmon{\mathrm{e}}
\def\comp{\mathrm{*}}
\newcommand{\ctx}{\mathrm{ctx}}
\newcommand{\sub}{\mathrm{sub}}
\newcommand{\ty}{\mathrm{ty}}
\newcommand{\tm}{\mathrm{tm}}
\def\nt{\mathrm{nat}}
\def\fun{\mathrm{fun}}

\title[Generalized Algebraic Theories and Categories with Families]{A Note on Generalized Algebraic Theories\\and Categories with Families}\author{Marc Bezem, Thierry Coquand, Peter Dybjer, Mart\'in Escard\'o}

\begin{document}

\maketitle

\begin{abstract}
We give a new syntax independent definition of the notion of a finitely presented generalized algebraic theory as an initial object in a category of categories with families (cwfs) with extra structure. To this end we define inductively how to build a valid signature $\Sigma$ for a generalized algebraic theory and the associated category $\Cwf_\Sigma$ of cwfs with a $\Sigma$-structure and cwf-morphisms that preserve $\Sigma$-structure on the nose.  Our definition refers to the purely semantic notions of {\em uniform family} of contexts, types, and terms. Furthermore, we show how to syntactically construct initial cwfs with $\Sigma$-structures. This result can be viewed as a generalization of Birkhoff's completeness theorem for equational logic. It is obtained by extending Castellan, Clairambault, and Dybjer's construction of an initial cwf. We provide examples of generalized algebraic theories for monoids, categories, categories with families, and categories with families with extra structure for some type formers of dependent type theory. The models of these are internal monoids, internal categories, and internal categories with families (with extra structure) in a category with families. Finally, we show how to extend our definition to some generalized algebraic theories that are not finitely presented, such as the theory of contextual categories with families.
\end{abstract}

\section{Introduction}

Martin-Löf type theory can be characterized in a syntax independent way as the initial category with families (cwf)  with extra structure for the type formers \cite{castellan:tlca2015,castellan:lmcs}. The main contribution of this note is a similar syntax independent characterization of the notion of finitely presented generalized algebraic theory as the initial cwf with extra structure.

Generalized algebraic theories (gats) were introduced by Cartmell in his PhD thesis \cite{cartmell:phd} as a dependently typed generalization of many sorted algebraic theories. Each gat is specified by a signature with (possibly infinite) sets of sort symbols, operator symbols, and equations. Cartmell's definition of gats \cite{cartmell:phd,cartmell:apal} is based on a notion of {\em derived rule} expressed in terms of a traditional syntactic system for dependent type theory. He also defines a notion of model whereby sort symbols are interpreted as families of sets.

Categories with families (cwfs) \cite{dybjer:torino} were introduced as a new notion of model of dependent type theory. Cwfs arise by reformulating the notion of category with attributes in Martin Hofmann's sense \cite{hofmann:csl}. The key point is that cwfs arise as models of a certain generalized algebraic theory closely related to Martin-Löf's substitution calculus \cite{martinlof:gbg92}. As such the notion of cwf becomes a useful intermediary between traditional syntactic systems for dependent type theory and a variety of categorical notions of model.

The gat of cwfs is thus a kind of idealized formal system of dependent type theory. In contrast to Martin-Löf's substitution calculus, and other syntactic systems for dependent type theory, it is {\em not} formulated in terms of grammars and inference rules for the forms of judgment of type theory. Instead it is formulated in terms of the sort symbols (corresponding to the judgment forms), operator symbols (corresponding to the formation, introduction, and elimination rules), and equations (corresponding to the equality rules for the type formers) of the gat. Some of the general reasoning (about equality, substitution, and assumptions) is taken care of by the underlying infrastructure of dependent types. This makes it possible to abstract away from details in the formulation of grammars and inference rules. In contrast to the various syntactic systems, the gat of cwfs has a canonical flavour. 

In this note we explore the interdependence between gats and cwfs. We already explained that cwfs can be defined as models of a gat. 
In the other direction, the notion of gat relies on the notion of cwf, in the sense that the latter models the underlying infrastructure of dependent types. 



\subsection*{Plan of the paper}

In Section 2 we recall the definition of the category $\Cwf$ of categories with families and morphisms preserving cwf-structure on the nose. Section 3 contains our main definition of a syntax independent notion of valid signature $\Sigma$ for a gat and the category $\Cwf_\Sigma$ of cwfs with a $\Sigma$-structure. In Section 4 we construct an initial object $\T_\Sigma$ in $\Cwf_\Sigma$. In Section 5 we show several examples of gats: for monoids, categories, cwfs, and cwfs with extra structure for one universe. We point out that cwfs with extra structure for gats of monoids, categories, cwfs are cwfs with an internal monoid, category, and cwf, respectively. We also sketch how to extend our approach to some countably presented gats, and show the example of contextual cwfs, a variant of Cartmell's contextual categories \cite{cartmell:phd,cartmell:apal}. Finally, in Section 6 we discuss related work, for example relating to Voevodsky's initiality conjecture \cite{voevodsky:initiality} and Altenkirch and Kaposi's quotient inductive-inductive types \cite{altenkirch:qiits}.

Our development can be formulated in a constructive set theory,
as described for instance by Aczel \cite{MR519801}, although the set theory
we use for formulating the notion of cwf with a $\Sigma$-structure is probably
much weaker. As emphasized by Voevodsky~\cite{voevodsky:initiality}, we study structures invariant
under {\em isomorphisms} and not under {\em equivalences}, and it is actually misleading
to call them ``category'' (and this is why Voevodsky used the term ``$C$-system''
for what Cartmell called ``contextual category'').
As he also noticed, this
important distinction between categories and notions invariant under isomorphisms becomes
precise in the setting of univalent foundations where not all collections of objects
are constructed from sets.


\subsection*{Remarks on terminology and notation}
Like Cartmell, we have chosen to use the term {\em sort symbol} from many-sorted universal algebra. However, in our semantic notion of signature sort symbols are interpreted as {\em type families} in a cwf. A cwf consists of a base category where the objects of the base category are (semantic) contexts and the morphisms are (semantic) substitutions. Moreover, we have a family-valued presheaf mapping contexts to families of (semantic) terms indexed by (semantic) types. Thus the reader should be aware of the mismatch between the word {\em sort} from universal algebra and the word {\em type} in the cwf semantics.

Another possible source of confusion is that cwfs appear on two different levels. In Section \ref{sec:def_cwf} we recall the definition of cwf in set-theoretic metalanguage, where we use $\Ty$ to denote the family of types indexed by contexts and $\Tm$ to denote the family of terms indexed by contexts and types. This notion of cwf is then used to define the semantic notions of signature and category of models of a gat. Then in Section \ref{gat-cwf} we define the gat of {\em internal cwfs}. This gat has sort symbols $\ty$ for {\em internal types} and $\tm$ for {\em internal terms} using lower case to highlight the difference from $\Ty$ and $\Tm$ in the model cwf. 

Furthermore, we often use the same notation both on the semantic and the syntactic level. For example, in Section \ref{gat-sig-mod},  where we are syntax independent, the letter $S$ denotes a semantic sort symbol, whereas in Section \ref{initial-gat}, where we construct the initial model, it denotes a syntactic sort symbol.

%
%
%

\section{Categories with families}\label{sec:def_cwf}

\subsection{The category of cwfs and strict cwf-morphisms}


\begin{definition}\label{def:catFam}
$\Fam$ is a category whose objects are
set-indexed families of sets, denoted as $(U_x)_{x\in X}$.
A morphism of $\Fam$ with source $(U_x)_{x\in X}$ and target $(V_y)_{y\in Y}$
consists of a re-indexing function $f: X\to Y$ together with a family
$(g_x)_{x\in X}$ of functions $g_x : U_x \to V_{f(x)}$. 
\end{definition}

The next step is to define the category $\Cwf$.
We split this definition in two: first the objects,
which are called \emph{categories with families}, in Definition~\ref{def:Cwfobj},
and then the morphisms in Definition~\ref{def:Cwfmor}.
Since $\Cwf$ has been developed as a categorical framework for the semantics of
type theory, much of the terminology (contexts, substitutions,
types, terms) refers to the syntax of type theory,
suggesting the intended interpretation of this syntax in the
so-called $\Cwf$-semantics.

The main novelty of this note is to use $\Cwf$ as a framework
to define a new notion of a generalized algebraic theory.
Contexts, substitutions, types, and terms also make
sense in relation to gats.

\begin{definition}\label{def:Cwfobj}
A category with families (cwf) consists of the following data:

\begin{itemize}
\item A category $\C$;

\item A $\Fam$-valued presheaf on $\C$, that is, a functor
$T : \Cop \to \Fam$;

\item A terminal object $1\in \C$, and unique maps
$\tuple{}_\Gamma \in \C(\Gamma, 1)$ for all objects $\Gamma$ of $\C$;

\item Operations ${\cext\,},~\tuple{\_,\_},~\p$ and $\q$
explained in the following paragraphs.
These four operations and their associated equations
are referred to as \emph{context comprehension}.
\end{itemize}

We let $\Gamma, \Delta,\ldots$ range over objects of $\C$,
and refer to them as \emph{contexts}.
We let $\delta, \gamma,\ldots$ range over morphisms,
and refer to them as \emph{substitutions}.
We refer to $1$ as the \emph{empty} context; the terminal maps
$\tuple{}_\Gamma$ represent the \emph{empty} substitutions.

If $T(\Gamma) = (U_x)_{x\in X}$, we write $\Ty(\Gamma)$ for the set $X$.
We call the elements of $\Ty(\Gamma)$ \emph{types in context $\Gamma$},
and let $A, B, C$ range over such types.
Furthermore, for $A \in \Ty(\Gamma)$, we write $\Tm(\Gamma, A)$ for the set $U_A$
and call the elements of $\Tm(\Gamma, A)$
\emph{terms of type $A$ in context $\Gamma$}.

For $\gamma : \Delta \to \Gamma$,
the functorial action of $T$ yields a morphism
\[
T(\gamma) \in  \Fam\left((\Tm(\Gamma, A))_{A\in \Ty(\Gamma)}, 
                (\Tm(\Delta, B))_{B\in \Ty(\Delta)}\right)
\]
consisting of a reindexing function $\_\,[\gamma] : \Ty(\Gamma) \to
\Ty(\Delta)$ referred to as \emph{substitution in types}, and for each $A\in
\Ty(\Gamma)$ a function $\_\,[\gamma] : \Tm(\Gamma, A) \to \Tm(\Delta,
A[\gamma])$, referred to as \emph{substitution in terms}.

Now we turn to the explanation of the operations
${\cext\,},~\tuple{\_,\_},~\p,~\q$.
Given $\Gamma \in \C$, $A \in \Ty(\Gamma)$, $\gamma : \Delta \to \Gamma$,
and $a\in \Tm(\Delta, A[\gamma])$, we have
\[
\Gamma \cext A \in \C
\quad\qquad
\p_{\Gamma, A} : \Gamma \cext A \to \Gamma
\quad\qquad
\q_{\Gamma, A} \in \Tm(\Gamma\cext A, A[\p_{\Gamma,A}])
\quad\qquad
\tuple{\gamma, a}_A : \Delta \to \Gamma \cext A.
\]
We call $\Gamma \cext A$ the \emph{extended} context
and $\tuple{\gamma, a}_A$ the \emph{extended} substitution.

The operations  ${\cext\,},~\tuple{\_,\_},~\p,~\q$
satisfy the following universal property:
$\tuple{\gamma, a}_A$ is the unique substitution satisfying
\[
\p_{\Gamma, A} \circ \tuple{\gamma, a}_A = \gamma
\qquad \text{and}\qquad
\q_{\Gamma, A} [\tuple{\gamma, a}_A] = a\,.
\]
We refer (colloquially) to $\p$ as the \emph{first projection},
and to $\q$ as the \emph{second projection}. 
{Note that the first equation implies that
$\Tm(\Delta,A[\p_{\Gamma,A}][\tuple{\gamma, a}]) = \Tm(\Delta,A[\gamma])$
so that $\q_{\Gamma, A} [\tuple{\gamma, a}]$ and $a$ are elements of the same set.}
Here and below, subscripts are omitted from ${\cext\,},~\tuple{\_,\_},~\p,~\q$
when they can be reconstructed from the context (no pun intended).
(End Definition~\ref{def:Cwfobj}.)
\end{definition}

A cwf is thus a structure $(\C,1,\tuple{},T,\cext\, , \tuple{\_,\_},\p, \q)$,
subject to equations, for the category and the presheaf, and universal
properties, formulated purely equationally, for the terminal object and for context comprehension.
The morphisms to be defined next preserve this structure,
even in a strict way, `on the nose'.
We often shorten the notation of a cwf to $(\C,T)$, or even just $\C$,
leaving the remaining structure implicit.

\begin{definition}\label{def:Cwfmor}
A \emph{(strict) cwf-morphism $F$ between cwfs $(\C,T_\C)$ and $(\D,T_\D)$}
consists of

\begin{itemize}

\item A functor $F_\fun : \C \to \D$;
\item A natural transformation $F_\nt : T_\C \Rightarrow (T_\D \circ F_\fun^\op)$;
\item The terminal object is preserved on the nose: $F_\fun(1_{\C}) = 1_{\D}$;
\item Context comprehension is preserved on the nose, see below.
\end{itemize}

Since $F_\nt$ is a natural transformation between $\Fam$-valued presheaves,
$F_\nt$ has a component for any object $\Gamma$ of $\C$, and
these components are morphisms in $\Fam(T_C(\Gamma),T_\D(F_\fun(\Gamma)))$.
Recall that morphisms in $\Fam$ consist of a reindexing function
and a family of functions. It is convenient to denote $F_\fun$,
all reindexing functions, as well as all members of the families of functions,
simply by $F$. Thus we have $F(A) \in \Ty_\D(F(\Gamma))$
and $F(a) \in \Tm_\D(F(\Gamma), F(A))$, for all $\Gamma$
and $A\in\Ty_\C(\Gamma)$ and $a\in \Tm_\C(\Gamma, A)$.

Naturality of $F_\nt$
amounts to preservation of substitution, {i.e.}, for all
$\gamma : \Delta \to \Gamma$ in $\C$, we have
\[
F(A[\gamma]) = F(A)[F(\gamma)] \qquad \qquad
F(a[\gamma]) = F(a)[F(\gamma)]\,.
\]

Last but not least, we turn to the preservation of context comprehension
on the nose, and require
\[
F(\Gamma\cext A) = F(\Gamma)\cext F(A) \qquad
F(\p_{\Gamma, A}) = \p_{F(\Gamma), F(A)} \qquad
F(\q_{\Gamma, A}) = \q_{F(\Gamma), F(A)}\,.
\]

Note that the universal property implies that
$F(\tuple{\gamma,a}) = \tuple{F(\gamma),F(a)}$.
The same is true for the terminal maps:
$F(\tuple{}_\Gamma) = \tuple{}_{F(\Gamma)}$.
(End Definition~\ref{def:Cwfmor}.)
\end{definition}

Small cwfs with strict cwfs-morphisms form a category, written $\Cwf$.

\section{Signatures and models of generalized algebraic theories}\label{gat-sig-mod}

We now come to the main point of this note.
We define how to build a valid gat signature $\Sigma$ and the associated
category $\Cwf_{\Sigma}$ of cwfs with a $\Sigma$-structure.
Each object of $\Cwf_{\Sigma}$ is a cwf with extra structure and
each morphism is a cwf-morphism preserving $\Sigma$-structure.
For this definition, we will need the following auxiliary notions.

A {\em uniform family of contexts} is a family $\Gamma = (\Gamma_{\C})$ with $\Gamma_\C$ a context in 
$\C$ for each $\C \in \Cwf_{\Sigma}$, such that
$F(\Gamma_\C) = \Gamma_\D$ for all morphisms $F \in \Cwf_{\Sigma}(\C,\D)$.
If $\Gamma$ is such a family, a {\em uniform family of types} over $\Gamma$ is a
family of types $A = (A_{\C})$ with $A_{\C}$ a type over $\Gamma_{\C}$ and
$F(A_{\C}) = A_{\D}$ for all morphisms $F \in \Cwf_{\Sigma}(\C,\D)$.
Finally, given $\Gamma$ and $A$, a {\em uniform family of terms} is a family
of terms $a = (a_{\C})$ with $a_\C \in \Tm_{\C}(\Gamma_{\C},A_{\C})$ such that
$F(a_{\C}) = a_{\D}$ for all morphisms $F \in \Cwf_{\Sigma}(\C,\D)$. 

\begin{remark}
Uniform families appear in Freyd's proof of the adjoint functor theorem \cite{freyd:abelian}, in Reynolds' \cite{reynolds:impredicative} and Reynolds and Plotkin's construction \cite{plotkin-reynolds} of an initial algebra for an endofunctor from an impredicative encoding of an inductive type, and in Awodey, Frey, and Speight's  \cite{awodey:impredicative} construction of an impredicative encoding of a higher inductive type. The common idea in these works is to first construct a weakly initial object and then the initial object is obtained by taking uniform families.
\end{remark}

\begin{definition}\label{def-sig-mod}
We define inductively (actually inductive-recursively) how to build a valid signature $\Sigma$ and the category $\Cwf_\Sigma$ of cwfs with a $\Sigma$-structure and cwf-morphisms that preserve $\Sigma$-structure. First, the base case:
\begin{description}
\item[The empty signature] The empty signature $\emptyset$ is valid and $\Cwf_\emptyset = \Cwf$.
\end{description}
Assume now that we have defined $\Sigma$ as a valid signature and the associated
category $\Cwf_{\Sigma}$.
Then we can add a new sort symbol, or a new operator symbol, or a new equation, to get a new valid signature,
as follows:
\begin{description}
\item[Adding a sort symbol]
  Let $\Gamma = (\Gamma_\C)$ be a uniform family of contexts indexed by $\C \in \Cwf_{\Sigma}$.
  Then we can extend $\Sigma$ with a new sort symbol $S$ relative to $\Gamma$, to obtain
  the gat $\Sigma' = (\Sigma,(\Gamma,S))$.
  The objects of $\Cwf_{\Sigma'}$ are pairs $(\C,S_{\C})$, where $\C$ is an object of $\Cwf_{\Sigma}$
  and $S_{\C} \in \Ty_\C(\Gamma_{\C})$.
  A morphism in $\Cwf_{\Sigma'}((\C,S_{\C}), (\D,S_{\D}))$
  is a morphism $F \in \Cwf_{\Sigma}(\C,\D)$ such that $F(S_\C) = S_\D$.
\item[Adding an operator symbol]
  If $\Gamma$ is a uniform family of contexts and $A$ a uniform family of
  types over $\Gamma$, 
  then we can extend $\Sigma$ with a new operator
  symbol $f$ relative to $\Gamma$ and $A$, to obtain
  the gat $\Sigma' = (\Sigma,(\Gamma,A,f))$.
  An object of $\Cwf_{\Sigma'}$
  is a pair $(\C,f_{\C})$ where $\C$ is an object in $\Cwf_{\Sigma}$ and $f_{\C} \in \Tm_\C(\Gamma_\C,A_\C)$.
  A morphism in $\Cwf_{\Sigma'}((\C,f_{\C}),(\D,f_{\D}))$ is a morphism $F \in \Cwf_{\Sigma}(\C,\D)$ such that $F(f_{\C}) = f_{\D}$
\item[Adding an equation]
  If $\Gamma$ is a uniform family of contexts,
  $A$ is a uniform family of types over $\Gamma$
  and $a,a'$ are uniform families of terms in $A$, 
  then we can extend $\Sigma$ with a new equation $a = a'$ relative to $\Gamma$ and $A$, to obtain
  the gat $\Sigma' = (\Sigma,(\Gamma,A,a,a'))$.
  In this case,
  $\Cwf_{\Sigma'}$ is a full subcategory of $\Cwf_{\Sigma}$. An object $\C$ in
  $\Cwf_{\Sigma'}$ is an object $\C$ of $\Cwf_\Sigma$ such that $a_{\C} = a'_{\C}$.
\end{description}

\end{definition}

This definition is {\em syntax independent}. In the next subsection we then show the purely {\em syntactic} construction of an {\em initial object} $\T_{\Sigma}$ in $\Cwf_{\Sigma}$ (for an arbitrary valid signature $\Sigma$) in terms of grammars and inference rules. A context in $\T_{\Sigma}$ will be an equivalence class $[ \Gamma ]$ of raw contexts, and similarly for substitutions, types, and terms. To give a uniform
  family of contexts $\Gamma_\C$ is then equivalent to giving a context $[ \Gamma ] \in \T_{\Sigma}$, since $\Gamma_\C = \inte{ [ \Gamma ] }_\C$ where $\inte{-}_\C$ is the interpretation morphism from $\T_\Sigma$ to $\C$.
Uniform families of types and terms arise from types and terms in $\T_\Sigma$ in a similar way.

We refer to Section \ref{monoids} where we show a simple example: the construction of a signature $\Sigma$ for internal monoids and its associated category of models $\Cwf_\Sigma$ of cwfs with an internal monoid. We also show how to construct the initial cwf with an internal monoid $\T_\Sigma$. 

\begin{remark}
Cartmell's notion of gat \cite{cartmell:phd,cartmell:apal} also makes it possible to stipulate equations between type expressions. However, neither of our examples makes use of this extra generality. In particular,  in Section \ref{sec:examples} we present the gat of cwfs with extra structure for $\N, \Pi$ and a first universe $\U_0$ without needing type equations. The reason is that
equations between types become equations between terms in our rendering
of dependent type theory as a gat. See Remark \ref{remark:typeequations} in Section \ref{sec:u-example} for more explanation.

Like Cartmell, we could consider gats with type equations, but we prefer not to make such equations part of our notion.

\end{remark}

\section{The construction of an initial object in $\Cwf_\Sigma$}\label{initial-gat}

In Section 3 we gave a {\em syntax independent specification} of a generalized algebraic theory as the initial object of the category $\Cwf_\Sigma$ of models of a (semantic) signature $\Sigma$. Now we show our main theorem: the {\em syntactic construction} of such an initial object $\T_\Sigma$. This construction is done in several steps. We first define the ``raw" syntactic expressions. Then we define four families of partial equivalence relations (pers) over those raw expressions, corresponding to the four equality judgments. The term model $\T_\Sigma$ is obtained by quotienting with these pers.

This theorem can be viewed as a generalization of Birkhoff's completeness theorem for equational logic \cite{birkhoff}:
\begin{theorem}
The category $\Cwf_\Sigma$ has an initial object $\T_\Sigma$,
for every valid signature $\Sigma$.
\end{theorem}

The construction of $\T_\Sigma$ will be by induction on the construction of $\Sigma$. It is based on construction of initial cwfs in \cite{castellan:tlca2015,castellan:lmcs} and we refer the reader to those papers for more details. Here we only provide a sketch and focus on how to extend the construction to $\T_\Sigma$.

For each $\Sigma$ we will define the following.
\begin{itemize}
\item
A grammar for the {\em raw syntax}, that is, raw contexts in $\RawCtx_\Sigma$, raw substitutions in $\RawSub_\Sigma$, raw types in $\RawTy_\Sigma$, and raw terms in $\RawTm_\Sigma$.
\item
A system of inference rules that generate four families of partial equivalence relations (pers) by a mutual inductive definition:
$$
\Gamma = \Gamma' \vdash_\Sigma
\qquad
\Gamma \vdash_\Sigma A = A'
\qquad
\Delta \vdash_\Sigma \gamma = \gamma' : \Gamma
\qquad
\Gamma \vdash_\Sigma a = a' : A
$$
where $\Gamma, \Gamma' \in \RawCtx_\Sigma, \gamma, \gamma' \in \RawSub_\Sigma, A, A' \in \RawTy_\Sigma,$ and $a,a' \in \RawTm_\Sigma$. These pers correspond to valid equality judgments of a variable-free version of dependent type theory with explicit substitutions based on the cwf-combinators. The ordinary judgments will be defined as the reflexive instances of these equality judgments. For example $\Gamma \vdash_\Sigma$ (meaning ``$\Gamma$ is a valid context") is defined as the reflexive instance $\Gamma = \Gamma \vdash_\Sigma$.
\item
A cwf $\T_\Sigma$ is then constructed from the equivalence classes of derivable judgments. For example, the contexts in $\T_\Sigma$ are equivalence classes $[\Gamma]$, such that $\Gamma \vdash_\Sigma$. We will show that $\T_\Sigma$ is a cwf with a $\Sigma$-structure, that is, an object of $\Cwf_\Sigma$.
\item
A $\Cwf_\Sigma$-morphism $\inte{-} : \T_\Sigma \to \C$ for every $\C \in \Cwf_\Sigma$. This is the {\em interpretation morphism}. This morphism is a partial function defined by induction on the raw syntax, that (whenever it is defined) maps raw contexts to contexts in $\C$, raw substitutions to substitutions in $\C$, raw types to types in $\C$, and raw terms to terms in~$\C$. We show that these partial functions preserve the partial equivalence relations so that we can define the interpretation morphism on the equivalence classes. Finally we show that it indeed is a $\Cwf_\Sigma$-morphism and the unique such into $\C$.
\end{itemize}

We begin with the construction for the base case: the {\bf empty signature} $\emptyset$.
\begin{itemize}
\item
We start with the {\em raw syntax}. This following grammar for raw contexts, raw substitutions, raw types, and raw terms.
\begin{eqnarray*}
\Gamma \in \RawCtx_\emptyset &::=& 1  \ |\ \Gamma\cext A\\
\gamma \in \RawSub_\emptyset  \ &::=& \gamma \circ \gamma \ |\ \id_\Gamma \ |\ \langle\rangle_\Gamma \ |\ \p_{A} \ |\ \langle \gamma, a \rangle_A\\
A \in \RawTy_\emptyset  &::=& A[\gamma]\\
a \in \RawTm_\emptyset  &::=& a[\gamma] \ |\ \qI_A
\end{eqnarray*}
These grammars generate a language of {\em cwf-combinators}.
\item
The system of inference rules is displayed in \cite{castellan:tlca2015,castellan:lmcs}. It is a system of {\em general rules}, rules for dependent type theory which come before we introduce any sort symbols and operator symbols and equations (or any rules for the type formers of intuitionistic type theory). We do not have room here to display them, but note that they can be divided into four groups:
\begin{itemize}
\item the per rules, amounting to symmetry and transitivity for the four forms of equality judgments;
\item preservation rules for judgments, amounting to substitution of equals for equals (an example of such a rule is the {\em type equality rule});
\item congruence rules for operators expressing that the cwf-combinators preserve equality;
\item conversion rules for the cwf-combinators.
\end{itemize}
\item
Note that the initial cwf $\T_\emptyset$ is rather uninteresting: its category of contexts contains only a terminal object (the empty context), and there are no types and terms. Nevertheless, the grammar and inference rules used in its definition form a starting point. The grammar for raw types and raw terms will be extended each time we add a new sort symbol or operator symbol, respectively. For each such new symbol and each new equation we will add a new inference rule. As a consequence we will generate a non-trivial $\T_\Sigma$.
\end{itemize}

Assume now for the induction step that we have defined the grammar, the inference rules, $\T_\Sigma$ and the interpretation morphism $\inte{-} : \T_{\Sigma} \to \C$ in $\Cwf_\Sigma$.
Let $\Sigma'$ be $\Sigma$ extended by a new sort symbol, a new operator symbol, or a new equation. We shall now explain how to define $\T_{\Sigma'}$.
\begin{description}
\item[Adding a sort symbol] If $\Gamma \vdash_\Sigma$, then we can introduce a new sort symbol $S$ in the context $\Gamma$ representing the sequence of types of the arguments of $S$.
\begin{itemize}
\item
We add a new production for raw types
$$
A ::= S
$$
to the productions for $\T_\Sigma$.
\item
We add  the inference rule
\begin{mathpar}
    \inferrule
    {}
    {\Gamma \vdash_{\Sigma'} S}
  \end{mathpar}
to the inference rules for $\T_\Sigma$.
\item
We define $S_{\T_{\Sigma'}} = [S]$, so that $ \T_{\Sigma'}$ has a $\Sigma'$-structure $(\T_\Sigma,[S])$.
\item
We extend the definition of the interpretation morphism $\inte{-}$  to an interpretation morphism $\inte{-}' : \T_{\Sigma'} \to \C$ by
$$
\inte{[S]}' = S_\C
$$
It follows  that this is a morphism in $\Cwf_{\Sigma'}$ and that it is unique.
\end{itemize}

\item[Adding an operator symbol] If $\Gamma \vdash_\Sigma A$, then we can introduce a new operator symbol $f$, where the context $\Gamma$ represents the sequence of types of the arguments and $A$ is the type of the result.
\begin{itemize}
\item
We add a new production for raw terms
$$
a ::= f
$$
to the productions for $\T_\Sigma$.
\item
We add  the inference rule
\begin{mathpar}
    \inferrule
    {}
    {\Gamma \vdash_{\Sigma'} f : A}
  \end{mathpar}
to the inference rules for $\T_\Sigma$.
\item
We define $f_{\T_{\Sigma'}} = [f]$, so that $ \T_{\Sigma'}$ has a $\Sigma'$-structure $(\T_\Sigma,[f])$.\item
We extend the definition of the interpretation morphism $\inte{-}$  to an interpretation morphism $\inte{-}' : \T_{\Sigma'} \to \C$ by
$$
\inte{[f]}' = f_\C
$$
It follows  that this is a morphism in $\Cwf_{\Sigma'}$ and that it is unique.
\end{itemize}

\item[Adding an equation] If $\Gamma \vdash_\Sigma a : A$ and $\Gamma \vdash_\Sigma a' : A$ we can introduce a new equation $a = a'$.
\begin{itemize}
\item
$\T_\Sigma'$ has the same productions as $\T_\Sigma$.
\item
We add  the inference rule
\begin{mathpar}
    \inferrule
    {}
    {\Gamma \vdash_{\Sigma'} a = a' : A}
  \end{mathpar}
to the inference rules for $\T_\Sigma$.
\item
$\T_{\Sigma'}$ is based on the same raw syntax as $\T_\Sigma$ but the equivalence relation has changed. To show that $\T_{\Sigma'} \in \Cwf_{\Sigma'}$ we just need to show that $[ a ] = [ a' ]$ but this follows from the inference rule $\Gamma \vdash_{\Sigma'} a = a' : A$.
\item
In order to define $\inte{-}'$ we first define the partial function on the raw syntax to be identical to the partial function on the raw syntax for $\inte{-}$. We then prove that this partial function preserves the extended partial equivalence relation and define $\inte{-}'$ on the new equivalence classes. It follows  that $\inte{-}'$ is unique.
\end{itemize}
\end{description}
This concludes the proof of the theorem. \qed

\section{Examples of generalized algebraic theories}\label{sec:examples}

We will now display the sort symbols, operator symbols, and equations for the generalized algebraic of internal monoids, internal categories and of internal cwfs. We will then show how to add operator symbols and equations when adding internal notions of $\Pi, \N$ and a universe closed under $\Pi$ and $\N$ to the gat of internal cwfs. The reason for prefacing these notions by the word ``internal" is that the models of the theories are internal monoids, categories, and cwfs in $\Cwf_\Sigma$ for the respective signatures for these theories. Moreover, internal monoids, categories, and cwfs in the cwf~$\Set$ are small monoids, categories, and cwfs, respectively. Note that the cwfs defined in Section 2 need not be small, and hence not internal cwfs in the cwf~$\Set$.

We begin by using the recipe in Definition \ref{def:Cwfmor} to construct the semantic signature for internal monoids and its associated category of models, that is, of cwfs with an internal monoid. We then follow the recipe in Definition \ref{initial-gat} and construct the initial cwf with an internal monoid. 

For ease of readability, we will only present the sort symbols, operator symbols, and equations in the remaining examples by using an informal notation with named variables, rather than the formal notations using cwf-combinators employed in Definitions \ref{def:Cwfmor} and \ref{initial-gat}.

Our final example is the gat of internal contextual cwfs, a variant of Cartmell's contextual categories. The contexts in such contextual cwfs come with a length $n$. We sketch how this can be axiomatized as a gat with countably many sort symbols $\ctx_n, \sub_n, \ty_n, \tm_n$ for an external natural number $n$ (and similarly for the operator symbols and equations). We also indicate how our framework can be extended to cover such gats.

\subsection{Internal monoids}\label{monoids}
 The one-sorted algebraic theory of monoids has two operator symbols,
$\idmon$ for identity and $\comp$ for composition, and associativity and identity laws as equations.
As any other one-sorted algebraic theory, the theory of monoids yields a
generalized algebraic theory. In ordinary notation with variables it might be rendered as follows, where $\Mon$ is the only sort:
\begin{eqnarray*}
&\vdash& \Mon\\
&\vdash& e : \Mon\\
x, y : \Mon &\vdash& \comp(x,y) : \Mon\\
y : \Mon &\vdash& \comp(\idmon,y) = y : \Mon\\
x : \Mon &\vdash& \comp(x,\idmon) = x : \Mon\\
x, y, z : \Mon &\vdash& \comp(\comp(x,y),z) = \comp(x,\comp(y,z)) : \Mon
\end{eqnarray*}

We now show how the corresponding semantic signature for monoids $\Sigma$ and its associated category of models $\Cwf_\Sigma$ are constructed step-wise in the sense of Definition \ref{def-sig-mod}.

As always, we begin with the empty signature $\emptyset$ and its category of models $\Cwf_\emptyset = \Cwf$.
\begin{description}
\item[Adding the sort symbol $\Mon$] Each cwf $\C$ has a chosen empty context (terminal object) $1_\C$. Since cwf-morphisms preserve empty contexts on the nose, $1 = (1_\C)$ is a uniform family of contexts in $\Cwf_\emptyset$. Hence we can introduce a new constant sort symbol $\Mon$ in the empty context. The resulting signature  is
$$
\Sigma_1 = (\emptyset, (1,\Mon))
$$
The objects of $\Cwf_{\Sigma_1}$  are pairs $(\C,\Mon_\C)$, where $\C$ is a cwf and $\Mon_\C \in \Ty_\C(1_\C)$.
\item[Adding the operator symbol for the identity]
Since, morphisms in $\Cwf_{\Sigma_1}$ preserve both empty contexts $1_\C$ and types $\Mon_\C$ on the nose, we have a uniform family of contexts $1 = (1_\C)$ and a uniform family of types $\Mon = (\Mon_\C)$ in $\Cwf_{\Sigma_1}$. Hence we can introduce a new constant operator symbol $e$ (the identity of the monoid).  The resulting signature  is
$$
\Sigma_2 = (\Sigma_1, (1,\Mon,e))
$$
The objects of $\Cwf_{\Sigma_2}$  are triples $(\C,\Mon_\C,e_\C)$, where $\C$ is a cwf, $\Mon_\C \in \Ty_\C(1_\C)$ and $e_\C \in \Tm_\C(1_\C,\Mon_\C)$.
\item[Adding the operator symbol for composition]
Again using that cwf-morphisms preserve all cwf-structure and $\Mon_\C$, we deduce that we have a uniform family of contexts $1.\Mon.\Mon[\p]$ and a uniform family of types $\Mon[\p][\p]$ in $\Cwf_{\Sigma_2}$. Thus we can add a binary operator symbol $\comp$. The resulting signature is
$$
\Sigma_3 = (\Sigma_2, (1.\Mon.\Mon[\p],\Mon[\p][\p],\comp))
$$
The objects of $\Cwf_{\Sigma_3}$  are quadruples $(\C,\Mon_\C,e_\C,*_\C)$, where $\C$ is a cwf, $\Mon_\C \in \Ty_\C(1_\C)$, $e_\C \in \Tm_\C(1_\C,\Mon_\C)$, and $*_\C \in \Tm_\C((1.\Mon.\Mon[\p])_\C,(\Mon[\p][\p])_\C)$.
\item[Adding the left identity law]
Furthermore, we extend the signature with the equations stating that $\idmon$ is a left identity as follows:
$$
\Sigma_4 = (\Sigma_3, (1.\Mon, \Mon[\p], \comp[\tuple{\tuple{\tuple{},\idmon[\tuple{}]},\q}], \q))
$$
The uniform family of context $1.\Mon$ expresses that the equation has one variable of type $\Mon$, the uniform family of types $\Mon[\p]$ expresses that the two sides of the equation have type $\Mon$, and the uniform families of terms $\comp[\tuple{\tuple{\tuple{},\idmon[\tuple{}]},\q}]$ and $\q$ express the two sides of the equation.
$\Cwf_{\Sigma_4}$ is the full subcategory of $\Cwf_{\Sigma_3}$ with objects $\C$
such that $(\comp[\tuple{\tuple{\tuple{},\idmon[\tuple{}]},\q}])_\C = \q_\C$.
\item[Adding the right identity and the associativity laws]
Finally we add the right identity equation and the associativity equation to get the signatures $\Sigma_5$ and $\Sigma_6$. We omit the details. 
\end{description}
We call the resulting signature for internal monoids $\Sigma = \Sigma_6$. The category $\Cwf_\Sigma$ is the category of cwfs with an {\em internal monoid}. This is a cwf-version of the notion of internal monoid which can be defined in any category with finite products. Ordinary (small) monoids come out as internal monoids in $\Set$, the cwf of small sets.

We also sketch the construction of the initial object $\T_\Sigma$ of $\Cwf_\Sigma$ following the recipe for introducing sort symbols, operators symbols, and equations in Section \ref{initial-gat}. (We omit the index $\Sigma$ in $\vdash_\Sigma$.)
\begin{description}
\item[Adding the sort symbol $\Mon$]
First, we have $1 \vdash$ for the empty signature, so we can
add a production for the constant sort symbol $\Mon$ and the inference rule:
\begin{eqnarray*}
1 &\vdash& \Mon
\end{eqnarray*}
For later use we infer $1.\Mon \vdash$ and, using $\p:1.\Mon \to 1$, $1.\Mon\vdash\Mon[\p]$,
so $1.\Mon.\Mon[\p] \vdash$.
\item[Adding the operator symbol for identity]
We then add a production for the constant operator symbol $\idmon$ and the inference rule:
\begin{eqnarray*}
1 &\vdash& \idmon : \Mon
\end{eqnarray*}
Again for later use we infer $1.\Mon\vdash \idmon[\p] :\Mon[\p]$.
Note that here $\p = \tuple{}$, the empty substitution $1.\Mon \to 1$,
since there is only one substitution $1.\Mon \to 1$.
\item[Adding the operator symbol for composition]
We then add a production for the binary operator symbol $\comp$.
Using another $\p: 1.\Mon.\Mon[\p] \to 1.\Mon$ (note the different type),
we can derive $1.\Mon.\Mon[\p] \vdash \Mon[\p][\p]$, so we can add the inference rule
\begin{eqnarray*}
1.\Mon.\Mon[\p] &\vdash& \comp : \Mon[\p][\p]
\end{eqnarray*}
Note that we project $\Mon$ on the right twice, reflecting that $\comp$ is binary.
\item[Adding the left identity law]
We can derive $1.\Mon \vdash \q : \Mon[\p]$. With some effort,
using previous inferences, we can  derive
$1.\Mon \vdash \comp[\tuple{\tuple{\tuple{},\idmon[\tuple{}]},\q}] : \Mon[\p]$.
Hence we can add the inference rule for the equation ($\idmon$ is a left identity):
\begin{eqnarray*}
1.\Mon &\vdash& \comp[\tuple{\tuple{\tuple{},\idmon[\tuple{}]},\q}] = \q : \Mon[\p]
\end{eqnarray*}
\item[Adding the right identity and the associativity laws]
We omit the details.
\end{description}
The resulting initial object $\T_\Sigma = \T_{\Sigma_6}$ is generated by a system of grammar and inference rules for dependent type theory with an internal monoid. In this theory we can prove statements such as 
$$
\Gamma \vdash a : \Mon
$$
stating that $a$ is a well-formed monoid expression in the context $\Gamma$ and
$$
\Gamma \vdash a = a': \Mon
$$
stating that $a = a'$ is a derivable equation between monoid expressions in the context $\Gamma$. Note that both contexts and monoid expressions use cwf-combinators and are variable-free. 

This example illustrates that gats indeed generalise ordinary algebraic theories.
The remaining examples will use dependent types in an essential way. However, for reasons of readability we will from now on only use ordinary notation with named variables. Hopefully, it is clear from the above how to formally construct the corresponding semantic signatures, categories of models, and initial models using cwf-combinators. For example, these constructions for the theory of internal categories are similar to the constructions for the theory of internal monoids.


\subsection{Internal categories} The gat of categories was one of Cartmell's motivating examples. It has the following sort symbols, operator symbols, and equations. Again, note that in our case the models are internal categories in a cwf. To emphasize the difference between the internal notions of category and cwf and the external notions (introduced in Section 2), our notation for sort symbols in the gat of internal cwfs use lower case letters ($\Obj, \Hom, \ty, \tm$). This is in contrast to the upper case letters for the external versions ($\Ty, \Tm$). We will however overload notation for operator symbols, and for example use $\circ$ both for the cwf-combinator and for the operator symbol in the gat of internal categories.

Sort symbols:
\begin{eqnarray*}
&\vdash& \Obj\\
\Delta, \Gamma : \Obj &\vdash& \Hom(\Delta,\Gamma)\\
\end{eqnarray*}

Operator symbols:
\begin{eqnarray*}
\Gamma : \Obj &\vdash& \id_\Gamma : \Hom(\Gamma,\Gamma)\\
\Xi,\Delta,\Gamma : \Obj, \gamma : \Hom(\Delta,\Gamma), \delta : \Hom(\Xi,\Delta) &\vdash&
\gamma \circ \delta : \Hom(\Xi,\Gamma)
\end{eqnarray*}

Equations:
\begin{eqnarray*}
\Delta, \Gamma : \Obj, \gamma : \Hom(\Delta,\Gamma) &\vdash& \id_\Gamma \circ \gamma = \gamma : \Hom(\Delta,\Gamma)\\
\Delta, \Gamma : \Obj, \gamma : \Hom(\Delta,\Gamma) &\vdash& \gamma \circ \id_\Delta = \gamma : \Hom(\Delta,\Gamma)\\
\Theta, \Xi,\Delta,\Gamma : \Obj, \gamma : \Hom(\Delta,\Gamma), \delta : \Hom(\Xi,\Delta), \xi : \Hom(\Theta,\Xi) &\vdash&
(\gamma \circ \delta) \circ \xi = \gamma \circ (\delta \circ \xi): \Hom(\Theta,\Gamma)
\end{eqnarray*}
Note that composition is officially an operator symbol with five arguments. In the official notation we should write $\gamma \circ_{\Xi,\Delta,\Gamma} \delta$, but we suppress the context arguments $\Xi,\Delta,\Gamma$. We will do so for some other operations too.

The rendering of the gat of categories in cwf-combinator language and the proof that it indeed yields a valid signature are similar to what they were for the gat of monoids. The inference rules for the two sort symbols in cwf-combinator language are
\begin{eqnarray*}
1 &\vdash& \Obj\\
1.\Obj.\Obj[\p] &\vdash& \Hom
\end{eqnarray*}
and the operator symbols for identity
\begin{eqnarray*}
1.\Obj &\vdash& \idmon : \Hom[\tuple{\id_{1.\Obj},\q_{1,\Obj}}]
\end{eqnarray*}
We omit the verbose cwf-renderings of the operator symbol for composition and the equations.

A cwf with extra structure for the generalized algebraic theory of categories is a cwf with an {\em internal category}. This is a cwf-based analogue of the usual notion of internal category in a category with finite limits. As shown by Martin Hofmann \cite{hofmann:csl,hofmann:cambridge}, every category with finite limits yields a category with attributes, and hence a cwf. However, not every cwf has finite limits. To achieve this we need more structure. As shown by Clairambault and Dybjer \cite{ClairambaultD11,ClairambaultD14} the 2-category of categories with finite limits is biequivalent to the 2-category of democratic cwfs that support $\Sigma$-types and extensional identity types.

An internal category in the cwf $\Set$ of small sets is a small category.


\subsection{Internal cwfs}\label{gat-cwf}

The gat of internal cwfs is obtained by extending the gat of internal categories with new sort symbols, operator symbols, and equations for a family valued functor, and then new operator symbols and equations for a terminal object, and context comprehension. We here rename the sort $\Obj$ of objects of the category of contexts to $\ctx$.

\subsubsection{The extension with a family valued functor}
\mbox{ }

Sort symbols:
\begin{eqnarray*}
\Gamma : \ctx &\vdash& \ty(\Gamma)\\
\Gamma : \ctx, A:\ty(\Gamma) &\vdash& \tm(\Gamma,A)
\end{eqnarray*}

Operator symbols:
\begin{eqnarray*}
\Gamma,\Delta : \ctx, A:\ty(\Gamma), \gamma : \Hom(\Delta,\Gamma) &\vdash&
A[\gamma] : \ty(\Delta)\\
\Gamma,\Delta : \ctx, A:\ty(\Gamma), \gamma : \Hom(\Delta,\Gamma), a:\tm(\Gamma,A) &\vdash&  a[\gamma] : \tm(\Delta,A[\gamma])
\end{eqnarray*}

Equations:
\begin{eqnarray*}
\Gamma : \ctx, A:\ty(\Gamma) &\vdash& A[\id_\Gamma] = A : \ty(\Gamma)\\
\Gamma : \ctx, A:\ty(\Gamma), a:\tm(\Gamma,A) &\vdash& a[\id_\Gamma] = a : \tm(\Gamma,A)\\
\Xi,\Delta,\Gamma : \ctx, \delta : \Hom(\Xi,\Delta), \gamma : \Hom(\Delta,\Gamma),
A:\ty(\Gamma) &\vdash& A[\gamma\circ\delta] = A[\gamma][\delta]: \ty(\Xi)\\
\Xi,\Delta,\Gamma : \ctx, \delta : \Hom(\Xi,\Delta), \gamma : \Hom(\Delta,\Gamma),
A:\ty(\Gamma), a:\tm(\Gamma,A) &\vdash&
a[\gamma\circ\delta] = a[\gamma][\delta]: \tm(\Xi,A[\gamma\circ\delta])
\end{eqnarray*}

\subsubsection{The extension with a terminal object}
No new sorts are required.

Operator symbols:
\begin{eqnarray*}
&\vdash& 1 : \ctx\\
\Gamma : \ctx &\vdash& \tuple{}_\Gamma : \Hom(\Gamma,1)
\end{eqnarray*}

Equations:
\begin{eqnarray*}
 &\vdash& \id_1 = \tuple{}_1 : \Hom(1,1)\\
\Gamma,\Delta : \ctx, \gamma : \Hom(\Delta,\Gamma) &\vdash&
\tuple{}_\Gamma\circ\gamma = \tuple{}_\Delta : \Hom(\Delta,1)
\end{eqnarray*}
(The latter two equations are better for term rewriting than the
obvious single one expressing the uniqueness of $\tuple{}_\Gamma$.)

\subsubsection{The extension with context comprehension}

No new sorts are required.

Operator symbols:
\begin{eqnarray*}
\Gamma : \ctx, A:\ty(\Gamma) &\vdash& \Gamma\cext A : \ctx\\
\Gamma,\Delta : \ctx, A:\ty(\Gamma), \gamma : \Hom(\Delta,\Gamma), a:\tm(\Delta,A[\gamma]) &\vdash& \tuple{\gamma,a} : \Hom(\Delta,\Gamma\cext A)\\
\Gamma : \ctx, A:\ty(\Gamma) &\vdash& \p: \Hom(\Gamma\cext A,\Gamma)\\
\Gamma : \ctx, A:\ty(\Gamma) &\vdash& \q: \tm(\Gamma\cext A,A[\p])
\end{eqnarray*}

Equations:
\begin{eqnarray*}
\Gamma,\Delta : \ctx, A:\ty(\Gamma), \gamma : \Hom(\Delta,\Gamma), a:\tm(\Delta,A[\gamma]) &\vdash& \p\circ\tuple{\gamma,a} = \gamma : \Hom(\Delta,\Gamma)\\
\Gamma,\Delta : \ctx, A:\ty(\Gamma), \gamma : \Hom(\Delta,\Gamma), a:\tm(\Delta,A[\gamma]) &\vdash& \q[\tuple{\gamma,a}] = a : \tm(\Delta,A[\gamma]) \\
\Gamma,\Delta,\Xi : \ctx, A:\ty(\Gamma), \gamma : \Hom(\Delta,\Gamma), a:\tm(\Delta,A[\gamma]), \delta : \Hom(\Xi,\Delta) &\vdash&
\tuple{\gamma,a} \circ \delta = \tuple{\gamma\circ\delta,a[\delta]} :
\Hom(\Xi,\Gamma\cext A) \\
\Gamma : \ctx, A:\ty(\Gamma) &\vdash&
\id_{\Gamma\cext A} = \tuple{\p,\q} : \Hom(\Gamma\cext A,\Gamma\cext A)
\end{eqnarray*}
(If $\p\circ\delta = \gamma$ and $\q[\delta]=a$, we get
$\tuple{\gamma,a}=\tuple{\p\circ\delta, \q[\delta]} = \tuple{\p,\q}\circ\delta =
\delta$, the uniqueness requirement of the universal property.
However, the equation for surjective pairing is not left-linear and with
a variable on one side, which is not good for rewriting.)

A cwf with extra structure supporting the generalized algebraic theory of cwfs is a cwf with an
\emph{internal cwf}. An internal cwf in the cwf $\Set$ of small sets is a {\em small cwf},
that is, it is a cwf in the ordinary sense (see Definition~\ref{def:Cwfobj})
except that it has small sets of objects, morphisms, types, and terms.

An example of a cwf with an internal cwf is provided by the cwf $\Set$ of small sets with an internal category of very small sets. We can make this precise if we work in set theory with two Grothendieck universes $\V_0 \in \V_1$. We call the members of $\V_1$ ``small sets" and the members of $\V_0$ ``very small sets". The category of contexts of the cwf $\Set$ is the usual category of small sets, by which we here mean that its objects are in $\V_1$. Moreover, the types are also the small sets in $\V_1$. To get an internal cwf, we interpret its sort of objects $\ctx$ as the small set $\V_0$ of very small sets, and the sorts of types $\ty(\Gamma)$ also as $\V_0$.


\subsection{Internal cwfs with $\Pi$-types}
We add three operator symbols in addition to the operator symbols for cwfs in Section 5.2 and 5.3:
\begin{eqnarray*}
\Gamma : \ctx, A : \ty(\Gamma), B : \ty(\Gamma.A)&\vdash& \Pi(A,B) : \ty(\Gamma)\\
\Gamma : \ctx, A : \ty(\Gamma), B : \ty(\Gamma.A), b : \tm(\Gamma.A, B) &\vdash& \lambda(b) : \tm(\Gamma,\Pi(A,B))\\
\Gamma : \ctx, A : \ty(\Gamma), B : \ty(\Gamma.A), c :  \tm(\Gamma,\Pi(A,B)), a : \tm(\Gamma, A) &\vdash& \app(c,a) : \tm(\Gamma, B[\tuple{\id,a}])
\end{eqnarray*}
(again omitting some of the official arguments)
and equations for $\beta, \eta$ (also omitting the context and type of the equality judgment)
 \begin{eqnarray*}
 \app(\lambda(b),a) &=& b[\tuple{\id,a}]\\
 \lambda(\app(c[\p],\q)) &=& c
\end{eqnarray*}
and commutation with respect to substitution.
\begin{eqnarray*}
\Pi(A,B)[ \gamma ] &=& \Pi(A [ \gamma ], B[ \gamma^+ ])\\
\lambda(b) [ \gamma ] &=& \lambda(b[\gamma^+ ])\\
\app(c,a) [ \gamma ] &=& \app(c[ \gamma ], a[ \gamma ] )
\end{eqnarray*}
where $\gamma^+ = \tuple{\gamma \circ \p, \q}$.

\subsection{Internal cwfs with $\Pi$ and $\N$}
Furthermore, we add the operator symbol
\begin{eqnarray*}
\Gamma : \ctx &\vdash& \N_\Gamma : \ty(\Gamma) 
\end{eqnarray*}
We also add operator symbols for $0, \s, \Rec$ and the equations for $\Rec$ and for commutativity with substitution, but omit the details. Note that the type of the primitive recursion operator $\Rec$ relies on the signature for $\Pi$-types.

\subsection{Internal cwfs with $\U_0$ closed under $\Pi$ and $\N$}\label{sec:u-example}
We add four more operator symbols
\begin{eqnarray*}
\Gamma : \ctx &\vdash& (\U_0)_\Gamma : \ty(\Gamma)\\
\Gamma : \ctx, a : \tm(\Gamma,(\U_0)_\Gamma) &\vdash& {\Ta_0}(a) : \ty(\Gamma)\\
\Gamma : \ctx &\vdash& (\N^0)_\Gamma : \tm(\Gamma,(\U_0)_\Gamma) \\
\Gamma : \ctx,
a : \tm(\Gamma,(\U_0)_\Gamma),
b :  \tm(\Gamma \cdot \Ta_0(a), (\U_0)_\Gamma))
&\vdash&
 \Pi^0(a,b) : \tm(\Gamma,(\U_0)_\Gamma)
\end{eqnarray*}
$(\U_0)_\Gamma$ is the universe (a type) relative to the context $\Gamma$; $\Ta_0$ is the decoding operation mapping a term in the universe to the corresponding type; $\N^0$ is the code for $\N$ in the universe, and $\Pi^0$ forms codes for $\Pi$-types in the
 universe. (Note that we have dropped the context argument of $\Ta_0$ and $\Pi^0$.)

We add the decoding equations:
\begin{eqnarray*}
\Ta_0(\N^0_\Gamma) &=& \N_\Gamma\\
\Ta_0(\Pi^0(a,b)) &=& \Pi(\Ta_0(a),\Ta_0(b))
\end{eqnarray*}
and the equations for preservation of substitution:
\begin{eqnarray*}
{(\U_0)}_\Gamma [ \gamma ] &=& {(\U_0)}_\Delta\\
\Ta_0(a) [ \gamma ] &=& \Ta_0(a[ \gamma ] )\\
\N^0_\Gamma [ \gamma ] &=&\N^0_\Delta\\
\Pi^0(a,b)[ \gamma ] &=& \Pi^0(a [ \gamma ], b[ \gamma^+ ])
\end{eqnarray*}

\begin{remark}\label{remark:typeequations}
Note that all equations are between {\em terms} in the gat of cwfs with extra structure for $\N, \Pi,$ and $\U_0$; we do not need the extra generality of stipulating type equations as discussed in the introduction. For example, $\Ta_0(\N^0_\Gamma) = \N_\Gamma$ is an equation between {\em internal} types, that is, terms of type $\ty(\Gamma)$.
\end{remark}

\begin{remark}
Also note that the gat for the universe is inevitably {\em \`a la Tarski} in the sense that we distinguish between types and terms in a cwf and we must have an operation decoding a term into a type. However, Martin-Löf's distinction between {\em \`a la Russell} and {\em \`a la Tarski} \cite{martinlof:padova} is a distinction between two different formulation of the raw syntax and inference rules of type theory.
\end{remark}

\subsection{A possible refinement to internal contextual cwfs}

Our treatment can be adapted to some non finitely presented gats.
If we have an increasing sequence of signatures $\Sigma_n$ we can consider their
union.
For instance, we can describe a gat of contextual cwfs \cite{castellan:lambek} (similar to Cartmell's contextual categories and Voevodsky's $C$-systems) by
the following stratification of the theory of cwfs. We replace the sort $\ctx$
by a sequence of sorts $\ctx_0,\,\ctx_1,\,\dots ,$ where $\ctx_n$ represents the sort
of contexts of length $n$ and a corresponding sequence of sorts
$\ty_n(\Gamma)$ for $\Gamma$ in $\ctx_n$
and $\tm_n(\Gamma,A)$ for $A$ in $\ty_n(\Gamma)$. Context extension $\Gamma.A$ is now in $\ctx_{n+1}$
if $A$ is in $\ty_n(\Gamma)$ and so on.
We also add {\em destructors}: we have
$\mathrm{ft}(\Gamma)$ in $\ctx_n$
and $\mathrm{st}(\Gamma)$ in $\ty_n(\mathrm{ft}(\Gamma))$
with $\Gamma = \mathrm{ft}(\Gamma).\mathrm{st}(\Gamma)$.
Similarly we have a stratification of the sort of substitutions
$\hom_{n,m}(\Delta,\Gamma)$ for $\Delta$ in $\ctx_n$ and $\Gamma$ in $\ctx_m$.
The resulting models are {\em internal contextual cwfs} in a cwf.

\begin{remark}
Generalized algebraic presentations of contextual categories (C-systems) have been given by Voevodsky \cite{voevodsky:c-systems} and Cartmell \cite{cartmell:gat-contextual}.
\end{remark}

\section{Related work}

Streicher \cite{streicher:semtt} defined {\em doctrines of constructions} (contextual categories with suitable extra structure) as a notion of model of the Calculus of Constructions. He also constructed a term model and remarked that it is an initial object in a category of doctrines of constructions.
Recently, Brunerie et al \cite{brunerie:initiality} presented a formalized proof in the Agda system that a formal system for Martin-Löf type theory forms an initial object in a category of contextual categories with extra structure for
the type formers.

More generally, Voevodsky \cite{voevodsky:initiality} outlined a new vision of the theory of syntax and semantics of dependent type theories. In this vision formal systems for dependent type theory are proved to be initial in suitable categories of models ({\em the initiality conjecture}). The above-mentioned contributions by Streicher and Brunerie et al are two examples of such characterizations. However, Voevodsky's aim was to go further and characterize a whole class of type theories and prove a general initiality result for them with the aim to form the basis for a general metatheory of dependent type theory. Our work can be viewed as a contribution to Voevodsky's programme, since we prove an initiality theorem for the whole class of finitely presented generalized algebraic theories. A characterization of a more general class of dependent type theories and their initial models has been proposed by Uemura \cite{uemura:general-framework}. Another related contribution is Palmgren and Vickers' \cite{palmgrenvickers} construction of initial models of essentially algebraic theories.

Altenkirch and Kaposi \cite{altenkirch:qiits} gave several examples of {\em quotient inductive-inductive types (qiits)}. Their main example is a definition of dependent type theory with $\Pi$-types and a universe, as a simultaneous definition in the Agda system \cite{agda-wiki} of the data types $\Ctx$ of contexts, $\Sub(\Delta,\Gamma)$ of substitutions, $\Ty(\Gamma)$, and $\Tm(\Gamma,A)$ of terms. Their definition is {\em inductive-inductive} \cite{nordvallforsberg:iids}, since the index sets of $\Sub, \Ty,$ and $\Tm$ are generated simultaneously, and as a consequence are not indexed inductive definitions in the usual sense where the index sets are fixed in advance. Furthermore, it is a quotient inductive-inductive type since they also have constructors for identity types, as in a {\em higher inductive type}.
There is a close relationship between this qiit and our initial internal cwf with $\Pi$-types and a universe. Like our definition, their qiit-definition uses cwf-combinators. Moreover, our sort symbols correspond to their formation rules (data type constructors), our operator symbols correspond to their introduction rules (constructors), and our equations correspond to their propositional identities. However, the fact that our equations are judgmental equalities while theirs are propositional identities is an important difference. As a consequence they need to use transport maps when moving between identical types. 

%

Kaposi, Kov{\'{a}}cs, and Altenkirch \cite{kaposi:qiits} developed a general theory of qiits. This includes a notion of signature for a qiit, the notion of an algebra of such a signature, and a construction of initial algebras. For these constructions they work in cwfs with $\Pi$-types, identity types, and a universe. This is in contrast to our work which is based on plain cwfs without extra structure for type formers. Although gats and qiits are related notions, neither is a generalization of the other. Gat is a basic notion independent of Martin-Löf type theory, whereas qiit is the latest in a series of generalizations of inductive type (inductive family, inductive-recursive type and family, inductive-inductive type, higher inductive type) extending intensional Martin-Löf type theory.



\subsection*{Acknowledgements} 
We are grateful to the anonymous referees for constructive criticism and pointers to related work. We would also like to thank Andrej Bauer and John Cartmell for further feedback.


\end{document}